\documentclass[12pt]{article}
\usepackage{graphicx}
\usepackage{amsmath}
\usepackage{amssymb}
\usepackage{amsthm}
\usepackage{latexsym}

\newcommand{\comment}[1]{}
\newtheorem{theorem}{Theorem}[section]
\newtheorem{lemma}[theorem]{Lemma}
\newtheorem{assumption}[theorem]{Assumption}
\newtheorem{proposition}[theorem]{Proposition}

\newtheorem{definition}[theorem]{Definition}

\newtheorem{remark}[theorem]{Remark}
\begin{document}

\date{}

\title{{\LARGE \textsf{Perfect simulation for unilateral fields }}}
\author{\textbf{Emilio De Santis, Mauro Piccioni}\\{\small
\textsl{Universit\`{a} di Roma \textit{La Sapienza\/},
Dipartimento di Matematica}}\\{\small \textsl{Piazzale Aldo Moro,
2 - 00185 Roma, Italia}} } \maketitle

\begin{abstract}
In this paper we consider two-point unilateral Markov fields on a
two-dimensional lattice as considered by Pickard \cite{Pi1},
Galbraith and Walley \cite{GW1,GW2}. We show that, under various
ergodicity conditions, they can be perfectly simulated in the
stationary state on any finite window. The techniques which are
used connect perfect simulation with oriented percolation through
suitable coupling constructions.
\end{abstract}

Mathematics Subject Classification 2000: Primary 60K35, 65C05;

 Secondary 60G60, 60G10.

Key words and phrases: Unilateral fields, Perfect simulation,
Oriented percolation.

\section{Introduction}

\label{Introduction}

In this paper we consider a particular class of random fields on
the two-dimensional lattice ${\mathbb N^+}^2$ and their stationary
extensions to ${\mathbb Z}^2$, the so-called {\it unilateral}
random fields. These fields appeared in the literature as models
for crystal growth \cite{WG73}; later they have been proposed also
for image analysis applications \cite {G}. Unilateral random
fields can be defined in general w.r.t. an oriented acyclic graph
structure
 (as in \cite{L}, where they are called Bayesian networks): for the sake of simplicity we
will consider the so-called two-point models on ${\mathbb Z}^2$,
which corresponds to the nearest neighbor structure where oriented
edges connect $(i,j)$ with $(i+1,j)$ and $(i,j+1)$. However the
reader will recognize that our basic idea works for more general
translation invariant graphs.

A unilateral two-point field model is constructed through a
transition kernel. We assume that the state space $E$ is Borel,
i.e. a Borel subset of a Polish space endowed with its Borel
$\sigma$-algebra. A (unilateral two-point) transition kernel $K$
on $E$ is a mapping $K:E^2 \rightarrow {\mathcal P}(E)$, where
${\mathcal P}(E)$ is the space of probability measures on $E$,
with the property that for any measurable set $B\subset E$ the
function $K(B|x,y)=K(x,y)(B)$ is measurable in the pair $(x,y) \in
E^2$. An $E$-valued unilateral field ${\mathbf X}^{\mathbf
x}=\{X_{i,j}^{\mathbf x} \in E, (i,j) \in {\mathbb N^+}^2 \}$ can
be constructed from $K$ for any choice of the boundary conditions
(b.c.'s) ${\mathbf x}=\{ x_{i,0}, x_{0,i}, i \in \mathbb N^+\}$.
The construction is accomplished by specifying consistently the
law of ${\mathbf X}_{\Lambda_{m,n}}^{\mathbf x}$ in any finite box
\begin{equation}
\Lambda_{m,n}=\{1,2,\ldots,m\}\times \{1,2,\ldots,n\}
\end{equation}
as being
\begin{equation}\label{field}
P({\mathbf X}_{\Lambda_{m,n}}^{\mathbf x}\in B)= \int_B
\prod_{i=1}^m\prod_{j=1}^n K(dx_{i,j}|x_{i-1,j},x_{i,j-1}),
\end{equation}
where $B$ is any measurable set in $E^{m\times n}$. A natural way
to simulate this field is to produce the variables
$X_{i,j}^{\mathbf x}$ in any sequential (total) order such that
each site $(i,j)$ comes after its {\it parents} $(i-1,j)$ and
$(i,j-1)$. A total order with this property will be called {\it
increasing} in the sequel. The reader will immediately notice the
similarity with discrete time Markov chains.

Next we randomize the b.c.'s, say with law $\mu$. We define the
law of the unilateral field ${\mathbf X}^{\mu}$ as
\begin{equation}
\label{randband}P({\mathbf X}_{\Lambda_{m,n}}^{\mu}\in A)=\int
P({\mathbf X}_{\Lambda_{m,n}}^{\mathbf x}\in A)\mu(d\mathbf x).
\end{equation}
We say that the b.c.'s are Markovian when the {\it horizontal
boundary} ${\bf X}_{\cdot , 0} = \{ X_{i,0} , i=1,2, \dots \}$ and
the {\it vertical boundary} ${\bf X}_{0, \cdot } = \{ X_{0, i} ,
i=1,2, \dots \}$ are Markov chains which are conditionally
independent given a common starting value $X_{0,0}$ (which can be
considered as the value of the field at the origin). For Markovian
b.c.'s the sequential simulation of the field ${\mathbf X}^{\mu}$
is easy, using any increasing order of the sites of ${\mathbb
N}^2$ (including boundary sites).

Pickard \cite {Pi1} studied the problem of determining laws $\mu$
such that the corresponding field ${\mathbf X}^{\mu}$ is
stationary, i.e. invariant under translations in ${\mathbb N}^2$.
Stationary unilateral fields can be extended to ${\mathbb Z}^2$ by
Kolmogorov's theorem. When $E$ is finite one can prove that
stationary unilateral fields can always be constructed. Using
again Kolmogorov's theorem, it suffices to take compatible
convergent subsequences of the sequence of averages
\begin{equation}\label{aver}
  \mu^L_{\Lambda_{m,n}}(\cdot )= \frac{1}{L^2} \sum_{t \in \Lambda_{L,L} }
  P( {\mathbf X^x}_{\Lambda_{m,n}+t} \in \cdot), \,\,\,\,\, L
  =1,2,\dots
\end{equation}
in any fixed finite box $\Lambda_{m,n}$, and use standard
compactness arguments.

For Markovian b.c.'s, Pickard established sufficient conditions
for stationarity of a unilateral field, which become also
necessary in the binary case. Pickard's conditions  are formulated
in terms of the joint distribution of the $4$-tuple of random
variables
$(X_{0,0}^{\mu},X_{1,0}^{\mu},X_{0,1}^{\mu},X_{1,1}^{\mu})$ (see
\cite {Pi1}), hence they involve $K$ and $\mu$. On the other hand,
Galbraith and Walley have proved that in the binary case any
positive two-point unilateral kernel $K$ has a unique Markovian
boundary law $\mu^K=\mu$ under which the columns
$(X_{0,0}^{\mu},X_{1,0}^{\mu})$ and
$(X_{0,1}^{\mu},X_{1,1}^{\mu})$ have the same law, and also the
rows $(X_{0,0}^{\mu},X_{0,1}^{\mu})$ and
$(X_{1,0}^{\mu},X_{1,1}^{\mu})$ have the same law  \cite {GW1}.
Since these two conditions are necessary for a stationary field,
 they could also reformulate Pickard's conditions
directly in terms of $K$.

However not all   stationary unilateral fields have Markovian
boundary laws, hence they are not necessarily easy to simulate
sequentially. The goal of this paper is to construct simulation
algorithms for the class of stationary unilateral fields which is
introduced next.

 We
say that a unilateral field $  {\mathbf X}^{\mathbf x}$  with
kernel $K$ is {\it ergodic} if there exists a stationary field
${\mathbf X'}=\{X_{i,j}' \in E, i,j=0,1,2,\dots \}$ such that for
every b.c.'s ${\bf x}$ and
 every pair $(m,n) \in {\mathbb N}^2$
\begin{equation}
\label{sta} \lim_{(h,l)\rightarrow \infty} ||{\mathcal L }({\bf
X}^{\bf x}_{(h,l)+\Lambda_{m,n}})-  {\mathcal L }({\bf
X}'_{\Lambda_{m,n}}) ||=0
\end{equation}
where we use ${\mathcal L }({\mathbf Y })$ to denote the law of a
random vector ${\mathbf Y } $, and  $||\cdot||$ is the total
variation. If the limit in (\ref{sta}) is uniform in the boundary
conditions $\mathbf x$ we say that the kernel is {\it uniformly
ergodic}. In the binary case Galbraith and Walley \cite {GW1} gave
a quite involved sufficient condition for the uniform ergodicity
of $K$.

We call the field ${\mathbf X'}$ appearing in (\ref{sta}) the {\it
equilibrium field} of the kernel $K$. It is easy to see that
${\mathbf X'}$ is a stationary unilateral field with kernel $K$,
and if it exists it is necessarily unique. Since its boundary law
$\mu'$ could be very hard to compute (even in the binary case), it
is not easy to simulate ${\mathbf X'}_{\Lambda}$ even on a small
box $\Lambda$. In fact, by (\ref {sta}) we know that we can only
{\it approach} ${\mathcal L }({\mathbf X'}_{\Lambda})$ by shifting
the box {\it sufficiently far away} from the boundary. The problem
is analogous to that of sampling a stationary Markov chain on a
finite window when the stationary distribution is not available
(as it happens in MCMC simulations); but in this case, in order to
simulate a box, we need to determine the whole joint equilibrium
law of its boundary.

Under a rather strong minorization condition on $K$ (Assumption
\ref{assume}) for a general Borel state space $E$, we show in the
next two sections how to implement a simulation algorithm which
produces a sample on any finite box, exactly distributed as the
equilibrium field. As a byproduct, we establish uniform ergodicity
of $K$. For one-dimensional discrete-time Markov processes the
first algorithm of this type was the CFTP algorithm of Propp and
Wilson \cite{PW}, see also \cite{Hag}. This work stimulated a wide
interest toward what is now called perfect simulation, see
\cite{Mad}. Our algorithm is based on an idea introduced by
Murdoch and Green \cite{MG} for discrete-time Markov processes. It
consists in coupling the whole family of unilateral fields, with
all possible boundary conditions, by using an underlying auxiliary
Bernoulli field. We show that the dependence from boundary values
is propagated only along increasing open paths in such a Bernoulli
field. In this way results from oriented percolation for
two-dimensional Bernoulli fields can be used in order to ensure
that the propagation stops with probability $1$ when the boundary
moves far away from the region to be simulated. In order to
simulate a square box of side $L$, the algorithm requires, in
addition to the $L^2$ variables of the box, an average of $O(L)$
additional random variables of the field. For other graphical
constructions used in the study of ergodicity see e.g.
\cite{Ferra}.

In Section \ref{varia}, in order to relax Assumption \ref{assume},
we study a more general class of algorithms working on blocks of
adjacent sites lying on selected diagonals. These algorithms work
under the more general Assumption \ref{assume3}, which is however
not equally easy to check. In Section \ref{aria}, using this block
algorithm, we present some examples of kernels for which
Assumption \ref{assume} fails but nevertheless a perfect
simulation algorithm can still be constructed.

\section{Coupling of unilateral fields}

\label{coupling}

In this section we introduce the auxiliary i.i.d. fields which
allow to couple, i.e. to represent on the same probability space,
the family of unilateral fields ${\mathbf X}^{\mathbf x}$ defined
in (\ref{field}), for all possible boundary conditions ${\mathbf
x} \in E^{\mathbb N^+}\times E^{\mathbb N^+}$. Moreover we
investigate when the values of these fields indicate that the
dependence from the boundary conditions of the field in some fixed
box is lost.

In the next two sections we suppose that the kernel $K$ satisfies
the following {\it minorization condition}.

\begin{assumption} \label{assume}
There exists a probability measure $\phi$ on $E$ and a positive
constant $\delta \geq \delta_0$ such that
\begin{equation}  \label{doblin}
K(A|y_1,y_2)\geq \delta \phi (A)
\end{equation}
for every measurable set $A$ and any pair $y_1,y_2 \in E$.
$\delta_0 $ is a positive constant which will be specified in the
next section.
\end{assumption}

If $\delta=1$ the field is trivially i.i.d. so we exclude this
possibility from now on.

\begin{remark}\label{finito} If $E$ is finite or countable, then provided
\begin{equation}  \label{sovra}
\tau(z) = inf \{ K( \{z \}| y_1, y_2 ) :\, y_1, y_2 \in E\}
\end{equation}
is not identically zero, the minorization condition (\ref{doblin})
is satisfied with
\begin{equation}  \label{inte}
\delta= \sum_{z \in E} \tau(z) > 0,
\end{equation}
and
\begin{equation}  \label{intr}
{\phi}(\{z\})=\frac{\tau(z)}{\delta}.
\end{equation}
\end{remark}

Under Assumption \ref{assume} it is immediately checked that
$$
H(\cdot | y_1,y_2)=\frac {K(\cdot|y_1,y_2)-\delta
\phi(\cdot)}{1-\delta}
$$
is a kernel on $E$. Since $E$ is Borel we can always define a
function $f:(0,1)\times E^2$ which is separately measurable in
each of its two arguments (see \cite{Ki}) with the property that
when $U$ is uniformly distributed in the interval $(0,1)$, then
$f({U; y_1,y_2})$ has the law $H(\cdot | y_1,y_2)$ for any pair
$(y_1,y_2) \in E^2$. Next we can prove the following

\begin{lemma}\label{labrulla}
Let $\{Z,V,U\}$ be mutually independent random variables with laws
\begin{equation}
P(Z=0)=1-P(Z=1)=\delta,\,\,\,\,\,\, V \sim \phi, \,\,\,\,\,\,U
\hbox{ uniform on } (0,1),\label{Z}
\end{equation}
and define
\begin{equation}
\label{recurs} g(z,v,u;y_1,y_2)=(1-z)v+zf(u;y_1,y_2).
\end{equation}
Then for any $y_1, y_2 \in E$
\begin{equation}\label{stoc}
g(Z,V,U;y_1,y_2) \sim K( \cdot |y_1,y_2).
\end{equation}
\end{lemma}
\begin{proof}\label{coup22}
For any $y_1 $, $y_2 \in E$ we can write the kernel as a mixture
\begin{equation}
\label{mixture} K (\cdot | y_1,y_2) = \delta \phi(\cdot)
+(1-\delta )H (\cdot | y_1,y_2)
\end{equation}
which is seen to be induced by the application of (\ref{recurs})
to $\{Z,V,U\}$.
\end{proof}

We say that the family of functions $\{g(\cdot;y_1,y_2),\,\, (y_1,
y_2) \in E^2\}$ realizes a {\it coupling} of the family of laws
$\{K(\cdot|y_1,y_2),\,\, (y_1, y_2) \in E^2\}$ on the probability
space where $Z$, $V$, and $U$ are defined. Notice that when $Z=0$
all the random variables $g(Z,V,U;y_1,y_2)$ take the same value
$V$, irrespectively of $(y_1, y_2) \in E^2$. In this case we say
that coupling occurs. When $E$ is finite or countable, with
$\delta$ and $\phi$ defined as in Remark \ref{finito}, the above
coupling is {\it maximal}, since the probability of coupling can
never exceed $\delta$. In fact, if $\{ X^{y_1,y_2}, \,\, (y_1,
y_2) \in E^2\}$ is any other coupling of $\{K(\cdot|y_1,y_2),\,\,
(y_1, y_2) \in E^2\}$, then
\begin{equation}
P(\exists x \in E: X^{y_1,y_2}=x,\,\, \forall (y_1,y_2) \in
E^2)\leq \sum_x \inf_{y_1,y_2} P(X^{y_1,y_2}=x)=\delta.
\end{equation}
Now let
$$\{\mathbf {Z},\mathbf {U},\mathbf {V}
\}=\{Z_{i,j},U_{i,j},V_{i,j} \}_{(i,j)\in {\mathbb Z}^2}$$ be
mutually independent i.i.d. fields with
$(Z_{i,j},U_{i,j},V_{i,j})\sim (Z,U,V)$ distributed as in
(\ref{Z}).

From the previous lemma we see that the probability space
supporting $\{\mathbf {Z},\mathbf {U},\mathbf {V} \}$ allows a
coupling of the fields ${\mathbf X}^{\mathbf x}$ for all boundary
conditions $\mathbf x \in E^{\mathbb N^+}\times E^{\mathbb N^+}$
by using the recursion
\begin{equation}
X_{i,j}^{\mathbf x}=g(Z_{i,j},V_{i,j},U_{i,j};X_{i-1,j}^{\mathbf
x},X_{i,j-1}^{\mathbf x}),\,\,\,\,\,\,\,\,\,i,j=1,2,\ldots
\label{rapp}
\end{equation}
along any increasing order of ${\mathbb N^+}^2$ starting from the
boundary conditions
\begin{equation}
X_{i,0}^{\mathbf x}=x_{i,0}\,,X_{0,i}^{\mathbf
x}=x_{0,i}\,\,\,\,\,\,\,\,\,\,\,i=1,2,\ldots.\label{rapp0}
\end{equation}

\medskip
The following definition will be useful in the sequel.

\begin{definition} \label{def1}
 For any
subset $S \subset {{\mathbb Z}^2}$ its {\emph { external
boundary}}, indicated with $ \overrightarrow{\partial} S$, is the
set of parents of some elements of $S$, which are not themselves
in $S$. Finally define the {\emph { internal boundary}} of $S$,
indicated with $\partial_I S$, as the set of vertices in $S$ with
at least one parent in $\overrightarrow{\partial} S$.
\end{definition}

Notice that the external boundary of the finite box
$\Lambda_{m,n}$ is the set of sites $\{(i,0),(0,j), i =1,\ldots,m,
j=1,\ldots,n \}$.

It is immediately seen that we can use the representation
(\ref{rapp}) to couple the field in any finite subset $S$ of
${\mathbb Z}^2$ for all possible boundary conditions prescribed on
its external boundary $\overrightarrow{\partial} S$, by recursion
along any increasing order of the sites in $S$. For any choice of
boundary conditions ${\mathbf x}=\{x_{i,j} \in E, (i,j) \in
\overrightarrow{\partial } S\}$ the unilateral field ${\bf
X}_S^{\bf x}$ on $S$, with b.c.'s ${\mathbf x}$ and transition
kernel $K$, can be therefore represented as
\begin{equation} \label{cxx}
{\bf X}_S^{\bf x}=G^S(\mathbf {Z}_{S},\mathbf {V}_{S},\mathbf
{U}_{S};{\mathbf x})
\end{equation}
where the function $G^S $ is suitably defined.

Furthermore we denote by $G^S_B(\mathbf {Z}_{S},\mathbf
{U}_{S},\mathbf {V}_{S};{\mathbf x})$ the projection of the
vector-valued function $G^S (\mathbf {Z}_{S},\mathbf
{U}_{S},\mathbf {V}_{S};{\mathbf x})$ on the sites belonging to
the subset $ B \subset S$.

By the direct inspection of (\ref{recurs}) we notice that if
$\mathbf {Z}_{\partial _{I}S}=\mathbf {0}$ then the recursion of
(\ref{rapp}) can be started with the configuration $\mathbf
{V}_{\partial _{I}S} $ irrespectively of the b.c.'s ${\mathbf x}
\in E^{\overrightarrow{\partial }S}$. We denote by
$\Gamma^S(\mathbf {Z}_{S},\mathbf {U}_{S},\mathbf {V}_{S}) $ the
resulting configuration of the field on $S$: hence
\begin{equation}\label{irresp}
{\mathbf {Z}}_{\partial _{I}S} = {\mathbf {0} }\Rightarrow
G^S({\mathbf {Z}}_{S},{\mathbf {U}}_{S},{\mathbf {V}}_{S};{\mathbf
x} )=\Gamma^S({\mathbf {Z}}_{S},{\mathbf {U}}_{S},{\mathbf
{V}}_{S}).
\end{equation}
The next lemma generalizes this situation.

\begin{lemma}
\label{R2} Let $S$ be a finite subset of ${\mathbb Z}^2$, and $B
\subset S$. If $\mathbf {Z}_{\partial _{I}B}=\mathbf {0}$, then
\begin{equation}\label{irresp2}
G^S_B(\mathbf {Z}_{S},\mathbf {U}_{S},\mathbf {V}_{S};{\mathbf x}
)=\Gamma^B(\mathbf {Z}_{B},\mathbf {U}_{B},\mathbf {V}_{B}),
\,\,\,\,\,\,\forall{\mathbf x} \in E^{\overrightarrow{\partial
}S}.
\end{equation}
\end{lemma}
\begin{proof} By the argument above the restriction of the field
on $B$ does not depend on its values on $\overrightarrow{\partial
}B$. By consequence it does not depend on the b.c.'s ${\mathbf x}$
on $\overrightarrow{\partial }S$.
\end{proof}

\section{Oriented percolation and coupling}

\label{oriented}

In this section the coupling made above is exploited to construct,
under Assumption~\ref{assume}, a perfect simulation algorithm for
the equilibrium distribution of an ergodic unilateral field in any
finite box $\Lambda_{m,n}$. Since we keep the box fixed, we write
it simply as $\Lambda$. Instead, the boundary conditions are
pushed far away from $\Lambda$, and the field is constructed using
a finite number of samples of a given realization of the auxiliary
fields $\mathbf {Z}$, $\mathbf {V}$ and $\mathbf {U}$. In
particular the auxiliary field $\mathbf {Z}$ will play a crucial
role, since it indicates when the construction becomes insensitive
to the boundary conditions.

We begin by defining an (increasing) {\it path} $\gamma$ joining
two vertices as a sequence of vertices $\{(i_k, j_k) \in {\mathbb
Z}^2\}_{k=0,1, \dots,m}$, with $m
>0$, such that for $k = 0,\dots,m-1$ either
\begin{equation}  \label{cam1}
i_{k+1} = i_k +1 , j_{k+1} = j_k \,\, \hbox{ or } \,\, i_{k+1} =
i_k , j_{k+1} = j_k +1.
\end{equation}
In this case we say that $(i_0, j_0)$ and $(i_m, j_m)$ are joined
by the path $\gamma$ of length $m$. Obviously two distinct
vertices $(a,b) \in {\mathbb Z}^2$ and $(c,d)\in {\mathbb Z}^2$
can be joined by a path if and only if $a\leq c $ and $b\leq d$.

Given a realization of the field $\mathbf {Z}=\{ Z_{i,j}\in
\{0,1\}, (i,j) \in {\mathbb Z}^2 \} $ we say that a path $\gamma
=\{(i_k, j_k ) \}_{k=0,1, \dots , m}$ is {\it open} (in the field
$\mathbf {Z}$) if $ Z_{i_k,j_k } =1 $ for $k=1, \dots , m$. By
convention we do not require that $Z_{i_0,j_0}=1$.

\begin{figure}[!ht]
\begin{center}
\includegraphics[width=7cm]{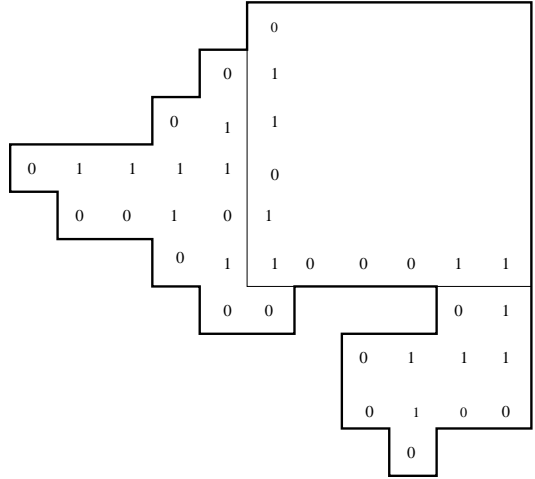}
\caption{a particular realization of $\omega (\Lambda)$.
}\label{fig0}
\end{center}
\end{figure}

We are now interested to study the following random subset of
${\mathbb Z}^2$
\begin{equation}\label{verticic}
\omega(\Lambda) = \omega(\Lambda, {\bf Z}) = \{ (i,j)\in {\mathbb
Z}^2: \exists \hbox{ an open path in $\mathbf {Z}$ joining } (i,j)
\hbox{ to } (k,l) \in
\partial_{I}\Lambda \} = \bigcup_{(k,l) \in \partial_I \Lambda} C_{k,l}
\end{equation}
where $ C_{k,l} $ is the set of vertices joined to $ (k,l) \in
\partial_I \Lambda$ by an open path in ${\mathbf Z}$ (see Figure
\ref{fig0} for a particular realization).

A fundamental result of oriented percolation for Bernoulli fields
(for a general reference on this subject see \cite{Du84}) ensures
that if the probability of $1$ at a site does not exceed a
critical value $0<p_c<1$ then $ C_{k,l}$ is finite with
probability $1$ (irrespectively of $(k,l)$ by translation
invariance). Moreover if such probability is strictly smaller than
$p_c$ then ${\mathbb E } | C_{k,l}| = {\mathbb E } | C_{0,0}| <
\infty $ (see again \cite{Du84}), which implies that
\begin{equation}\label{box}
{\mathbb E }| \omega ( \Lambda )|  \leq  |\partial_I \Lambda |
\cdot  {\mathbb E } | C_{0,0}|.
\end{equation}
Since an estimate from below is trivially obtained using the bound
\begin{equation}\label{below}
  |\omega (\Lambda)| \geq \sum_{(i,j ) \in \partial_I \Lambda }
  Z_{i,j},
\end{equation}
we conclude in this case that $ {\mathbb E }| \omega ( \Lambda )|
= O( |\partial_I \Lambda| ) $.

{\bf From now on we set $\delta_0=1-p_c$ in
Assumption~\ref{assume}.} Thus
\begin{proposition}
\label{percola} Under Assumption~\ref{assume} the random subset
$\omega(\Lambda)$ is finite almost surely. Moreover if $\delta >
\delta_0$ then the mean value of $|\omega(\Lambda)|$ is $O(
|\partial_I \Lambda | )$.
\end{proposition}

Since we are interested in bounding $\delta_0$ from above, we need
to bound $p_c$ from below. For example in \cite{Gra} it is proved
that $p_c\geq 0.682 \dots $, thus for $\delta \geq 0.317 \dots $
Assumption~\ref{assume} is verified.

Let us define the random set $B(\Lambda)= \Lambda \cup \omega
(\Lambda)$. Then we can prove the following
\begin{lemma}\label{uscita}
For any $(i,j) \in \partial_I B(\Lambda)$ it is $Z_{i,j}=0$.
\end{lemma}
\begin{proof}
We first notice that the internal boundary $\partial_I B(\Lambda)$
is the disjoint union of $\{(i,j): (i,j) \in
\partial_I \Lambda, Z_{i,j}=0\}$ and $\partial_I \omega
(\Lambda)$. Therefore, to prove the lemma, we need only to show
that $Z_{i,j}=0$ for any $(i,j) \in
\partial_I \omega(\Lambda)$. This is due to
the fact that if $Z_{i,j}=1$ then both the parents $(i-1,j)$ and
$(i,j-1)$ are in $\omega (\Lambda)$ since the existing open path
from $(i,j)$ to a site in $\partial_I \Lambda$ can be extended in
both directions. But this is absurd, since $(i,j) \in
\partial_I \omega(\Lambda)$.
\end{proof}

At this point the main result of the paper can be proved. As
before the function $\Gamma^B_A$ denotes the restriction of
$\Gamma^B$ to the sites in $A \subset B$.

\begin{theorem}\label{main}
Under Assumption~\ref{assume} the kernel $K$ is uniformly ergodic.
Moreover for any finite box $\Lambda=\Lambda_{m,n}$
\begin{equation}\label{boh}
\Gamma^{B(\Lambda)}_{\Lambda}(\mathbf {Z}_{B(\Lambda)},\mathbf
{V}_{B(\Lambda)},\mathbf {U}_{B(\Lambda)})\sim {\bf X'}_{\Lambda},
\end{equation}
where ${\bf X'}$ is the equilibrium field with kernel $K$.
\end{theorem}
\begin{proof}

By Proposition \ref{percola} $B=B(\Lambda)$ is almost surely
finite. Next, for any finite $S$ such that $\Lambda \subset S$,
let $Q^{\Lambda}_S$ be the event $\{B(\Lambda) \subset S\}$. From
Lemma \ref{R2} and Lemma \ref{uscita} we know that if
$Q^{\Lambda}_S$ occurs, for any $\mathbf x \in
E^{\overrightarrow{\partial }S}$
\begin{equation}
\label{funz} \mathbf X^{\mathbf x}_{S
|{\Lambda}}=G^S_{\Lambda}(\mathbf {Z}_{S},\mathbf {V}_{S},\mathbf
{U}_{S}; \mathbf x) = \Gamma^B_{\Lambda}(\mathbf {Z}_{B},\mathbf
{V}_{B},\mathbf {U}_{B}).
\end{equation}

Next, let $\{S_k \supset \Lambda, k=1,2,\dots\}$ be any sequence
of finite subsets increasing to the set of all sites $(h,l)\in
{\mathbb Z}^2$ such that $h\leq i$ and $j \leq l$. Then, under
Assumption~\ref{assume}, the sequence $\{Q^{\Lambda}_{S_k}\}$
increases to an event of probability $1$. This means that with
probability $1$ $B(\Lambda) \subset S_k$ eventually with $k$,
hence
\begin{equation}
\label{backward} G^{S_k}_{\Lambda}(\mathbf {Z}_{S_k},\mathbf
{V}_{S_k},\mathbf {U}_{S_k}; \mathbf
x_k)=\Gamma^B_{\Lambda}(\mathbf {Z}_{B},\mathbf {V}_{B},\mathbf
{U}_{B}) \,\,\,\, \forall \mathbf x_k \in
E^{\overrightarrow{\partial }S_k}
\end{equation}
eventually in $k$.

Next let us consider any  sequence $\{( h_k, v_k) k=1,\ldots\}$
converging towards $\infty$ and construct the sequence of boxes
$$
C_k :=\{1,\ldots,h_k+m\}\times \{1,\ldots,v_k+n\} \subset ( h_k,
v_k) +\Lambda = : \Lambda_k  .
$$
It is clear that for any b.c.  ${\mathbf x}_k \in
E^{\overrightarrow{\partial }C_k}$
\begin{equation}\label{legeq}
X^{{\mathbf x}_k}_{(h_k,v_k)+\Lambda}=G^{C_k}_{\Lambda_k}(\mathbf
{Z}_{C_k},\mathbf {V}_{C_k},\mathbf {U}_{C_k}; {\mathbf x}_k) \sim
G^{S_k}_{\Lambda}(\mathbf {Z}_{S_k},\mathbf {V}_{S_k},\mathbf
{U}_{S_k}; {\mathbf x}_k)
\end{equation}
where
$$
S_k=\{-h_k+1,\ldots,m\}\times \{-v_k+1,\ldots,n\},
$$
with the b.c.'s ${\mathbf x}_k$ at the r.h.s. of (\ref{legeq})
defined on
$$
{\overrightarrow{\partial }S_k} =
\{(-h_k,-v_k+1),\ldots,(-h_k,m)\}\bigcup
\{(-h_k+1,-v_k),\ldots,(m,-v_k)\}.
$$
The r.h.s of (\ref{legeq}) is equal to $\Gamma^B_{\Lambda}(\mathbf
{Z}_{B},\mathbf {V}_{B},\mathbf {U}_{B})$ on $Q^{\Lambda}_{S_k}$.
By consequence, denoting by $I_A $ the indicator of $A$
\begin{equation}\label{barc}
||{\mathcal L }( X^{{\mathbf x}_k}_{(h_k,v_k)+\Lambda} ) -
{\mathcal L }( \Gamma^{B}_{\Lambda} ) ||=||{\mathcal L }(
G^{S_k}_{\Lambda} ) - {\mathcal L }( \Gamma^{B}_{\Lambda} ) ||=
\sup_A {\mathbb E } ( | I_A ( G^{S_k}_{\Lambda} ) -   I_A (
\Gamma^{B}_{\Lambda} )| ) \leq 2P({Q_{S_k}^{\Lambda}}^c),
\end{equation}
from which we deduce uniform convergence in variation because $
\lim_{k \to \infty } P({Q_{S_k}^{\Lambda}}^c) =0 $. Since, by
varying $\Lambda$, this specifies a stationary and compatible
family, by Kolmogorov's Theorem it is a realization of a
stationary field ${\mathbf X'}$ in the finite box $\Lambda$.
\end{proof}

To summarize we have constructed a sampling scheme for the
equilibrium field ${\mathbf X'}$ on a finite box $\Lambda $, based
on a random but a.s. finite number of samples from the fields
$\mathbf {Z},\mathbf {V},\mathbf {U}$. In fact, by Proposition
\ref{percola} this number is proportional to $|\Lambda|$ plus $O(
|\partial_I \Lambda | )$ in the average.

We conclude the section by discussing a possible implementation of
the algorithm. We construct $\omega(\Lambda)$ by backward
induction in the following way. Let
$$
\Delta_0=\{(i,j):(i,j) \in \partial_I \Lambda, Z_{i,j}=1\}
$$
Then for $ k \geq 1$ we determine $\Delta_{k+1}$ from $\Delta_k$
as
$$
\Delta_{k+1}=\{(i,j): \{(i+1,j),(i,j+1)\} \cap \Delta_k \neq
\emptyset , Z_{i,j}=1\}
$$
until the index $k_{max}$ such that $\Delta_{k_{max}+1}=\emptyset$
for the first time, which is finite with probability $1$ if
$\delta\leq \delta_0$. Then $\omega (\Lambda)$ is the union of
$\cup_{k=1}^{k_{max}}\Delta_k$ and its external boundary. Then
order the sites in $\omega (\Lambda)$ by starting with those in
such an external boundary, then those in $\Delta_{k_{max}}$, next
those in $\Delta_{k_{max}-1}$ which were not already in
$\Delta_{k_{max}}$, and so on.  After having totally ordered all
sites in $\omega(\Lambda)$, we start with those in $\Lambda$.  No
matter which order is chosen within each of these regions, the
recursion along this total order will allow to compute
${\Gamma}^B_{\Lambda}(\mathbf {Z}_{B},\mathbf {V}_{B},\mathbf
{U}_{B})$.

\section{A class of more general block algorithms} \label{varia}

In the previous section we have shown how to construct a perfect
simulation scheme for some stationary ergodic unilateral fields on
${\mathbb Z}^2$. We have used percolation arguments on a suitable
auxiliary Bernoulli field to show that the algorithm works under
the "sufficiently large" minorization condition Assumption
\ref{assume}. On the other hand, in the one-dimensional case the
same idea leads to an algorithm for the exact simulation of a
sample from a stationary discrete-time Markov chain that works
under any non trivial minorization condition, since in any non
degenerate one-dimensional Bernoulli field clusters are always
finite.

More generally, the Multigamma coupler of Murdoch and Green
\cite{MG} extends the above idea to cover the whole class of
uniformly ergodic kernels by considering a suitable power of the
kernel. In fact, by Theorem 16.0.2 in \cite{MT} uniformly ergodic
kernels are characterized by a minorization condition on some
power $K^m$ of the kernel $K$, which means that we can apply the
same algorithm to the $m$-skeleton chain $\{X_{km}, k=0,1,\ldots
\}$. Moreover  provided $m$ is large enough the value of $\delta$
in the minorization condition can be taken arbitrarily close to
$1$. More generally, it has been proved in \cite{FT} that a
vertical backward coupling time exists only for such a class of
Markov chains.

Based on these observations, in this section we try to extend the
previous results  to a wider class of unilateral fields by
considering suitable skeleton fields, which in general consist of
{\it blocks} of sites.

By a diagonal of ${\mathbb Z}^2$ we mean a set
$$
D_h=\{(i,j)\in {\mathbb Z}^2:i+j=h\}
$$
for $h \in {\mathbb Z}$. The distance between $D_h$ and $D_k$ is
defined as $|h-k|$. The binary field used to indicate the region
where the simulation have to be performed is in general defined
over a new lattice associated to blocks of $l$ adjacent sites
lying on diagonals at distance $(d-1)l$ one from the other, for
given integers $l=1,2,\ldots$ and $d=2,3,\ldots$. We start by
defining the block
\begin{equation}\label{siti}
  B_{0,0}=\{(1,l),(2,l-1)\ldots,(l,1)\}
\end{equation}
and, for any pair $(i,j) \in {\mathbf Z}^2$ such that $i+j=0$
$mod(d-1)$, define its translates
\begin{equation}\label{otis}
  B_{i,j}=B_{0,0}+(il,jl)
\end{equation}

\begin{figure}[!ht]
\begin{center}
\includegraphics[width=7cm]{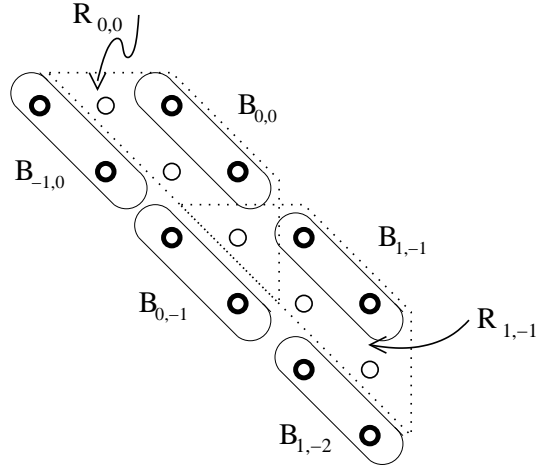}
\caption{Blocks with $l=2$ and $d=2$.}  \label{fig1}
\end{center}
\end{figure}

These blocks are taken to be the vertices of a new graph
${\mathcal G}_d=({\mathcal V}_d, {\mathcal E}_d)$. An oriented
edge connects $B_{h,k}$ with $B_{h+i,k+d-1-i}$ for
$i=0,1,\ldots,d-1$, for any pair $(h,k) \in {\mathbb Z}^2$ such
that $h+k=0$ $mod(d-1)$. In the following we will directly
identify ${\mathcal V}_d$ with such a subset of sites of ${\mathbb
Z}^2$ and refer to the sites $(h+i,k-i-d+1)\in {\mathcal V}_d$ as
the $d$ parents of the site $(h,k) \in {\mathcal V}_d$, for
$i=0,1,\ldots,d-1$. It is clear that for $l=1, d=2$ we get that
the graph ${\mathcal G}_d$ is the original lattice ${\mathbb
Z}^2$; more generally, for larger values of $l$ and $d=2$, the
graph ${\mathcal G}_d$ is isomorphic to ${\mathbb Z}^2$. In the
next section we will devote Example 2 to show that it may be
convenient to let $l$ grow, rather than $d$.

Next let ${\cal S}=\{s_m, m=0,1,\ldots,(d-1)l \}$ be a family of
coupling functions, i.e. for $m=0,1,\ldots,(d-1)l$,
$s_m:(0,1)\times E^2 \rightarrow E$ is separately measurable in
each of its arguments and such that
\begin{equation}\label{couple}
  s_m(U;x_1,x_2)\sim K(\cdot|x_1, x_2),
\end{equation}
when $U$ is a uniformly distributed random variable taking values
in $(0,1)$. Since $E$ is Borel the function $g$ defined in
(\ref{recurs}) is a particular example of (\ref{couple}), since
all the auxiliary random variables $Z$, $V$ and $U$ which are
required can be constructed as functions of a single uniformly
distributed random variable. The construction of the field
${\mathbf X}^{\mathbf x}_S$ with kernel $K$ over any finite region
$S$ with b.c. $\mathbf x$ on $E^{\overrightarrow{\partial} S}$ is
then performed by using a vector ${\mathbf U}_S$ of independent
uniformly distributed random variables, iterating along an
increasing total order of the sites in $S$ the recursion
\begin{equation}\label{recu}
  X_{i,j}^{\mathbf x}=
  s_{m(i,j)}(U_{i,j};X_{i-1,j}^{\mathbf x},X_{i,j-1}^{\mathbf x})
\end{equation}
where
$$
m(i,j)=i+j-l-1,\,\,\,mod \{(d-1)l\},
$$
starting from
$$
X^{\mathbf x}_{\overrightarrow{\partial} S}=\mathbf x.
$$
Notice that the coupling function $s_m$ used to construct the
value of the field at a given site $(i,j)$ is allowed to depend on
$m(i,j)$, the distance of the diagonal where $(i,j)$ lies from the
``previous" diagonal of blocks. In Example 1 we will show the
usefulness of allowing couplings depending on the diagonal.

Next, for any fixed pair of integers $(h,k) \in {\mathcal V}_d$,
consider the trapezoidal region
$$
R_{h,k}=\{(i,j): i\leq l(h+1), j \leq l(k+1), (h+k-d+2)l+2\leq i+j
\leq (h+k+1)l+1 \}.
$$
The reason for defining this region is that $B_{h,k} \subset
R_{h,k}$ and
$$
\overrightarrow{\partial} R_{h,k}=\cup_{i=0}^{d-1}B_{h-i,k+i-d+1}.
$$
Therefore we may represent the field ${\mathbf X}^{\mathbf x}_
{R_{h,k}}$ as a function of the vector ${\mathbf U}_{R_{h,k}}$
with i.i.d. components and the boundary values on the parent
blocks
$$
\mathbf x=({\mathbf x}_{B_{h-i,k+i-d+1}}, i=0,\ldots,d-1)
$$
in particular
\begin{equation}\label{mah}
  \mathbf X^{\mathbf x}_{B_{h,k}}=F^{\cal S}(\mathbf U_{R_{h,k}};{\mathbf x}).
\end{equation}
where $F^{\cal S}$ is defined through the recursive application of
(\ref{recu}) along an increasing total order of the sites in
$R_{h,k}$.

We are now ready to make the following general assumption.
\begin{assumption}\label{assume3}
There exists a binary field
$$
W_{h,k}=\psi (\mathbf U_{R_{h,k}}), \,\,\,\,\,(h,k) \in {\mathcal
V}_d
$$
with
\begin{equation}\label{jam} P(W_{h,k}=0)=\tilde \delta >
{{d-1}\over{d}},
\end{equation}
such that $\{W_{h,k}=0\}$ implies
\begin{equation}\label{ham}
F^{\cal S}(\mathbf U_{R_{h,k}};{\mathbf x})=\Phi (\mathbf
U_{R_{h,k}}),\, \forall {\mathbf x} \in
E^{\overrightarrow{\partial} R_{h,k}}
\end{equation}
for some measurable function $\Phi$.
\end{assumption}

We kept Assumption \ref{assume3} quite general in order to
accomodate various possible definition of the field $\mathbf
W=\{W_{h,k}, (h,k) \in {\mathcal V}_d\}$, for a given family $\cal
S$ of coupling functions.

In principle (at least when the state space $E$ is countable) the
field $\mathbb W$ can be directly defined to have the value zero
if and only if $F^{\cal S}(\mathbf U_{R_{h,k}};{\mathbf x})$ does
not depend on ${\mathbf x}\in E^{\overrightarrow{\partial}
R_{h,k}}$, and in this case $\Phi (\mathbf U_{R_{h,k}})$ is equal
to such a common value. This means that in order to conclude that
$W_{h,k}=0$ we have to check that the realizations of the field
$\mathbf X^{\mathbf x}_{B_{h,k}}$ started from all the possible
b.c.'s ${\mathbf x} \in E^{\overrightarrow{\partial} R_{h,k}}$
collapse into a single value. Thus, within a single region
$R_{h,k}$ this is similar to the original Propp and Wilson
coupling from the past algorithm \cite{PW}.

A computationally less demanding choice is to define recursively
the random subset of $E$
\begin{equation}\label{range}
I_{i,j}=s_{m(i,j)}(U_{i,j};I_{i,j-1}\times I_{i-1,j}),
\end{equation}
for any site in $(i,j) \in R_{h,k}$, starting from $I_{i,j}=E$ for
$(i,j) \in  \overrightarrow {\partial} R_{h,k}$. Finally we define
$W_{h,k}=0$ if $I_{i,j}$ is a singleton for all $(i,j) \in
B_{h,k}$, and $\Phi ({\mathbf U}_{R_{h,k}})$ is then equal to its
unique element. Another option for the definition of $\mathbf W$
will be presented in Example 1.

A generalization of the coupling (\ref{recu}) is possible by
allowing $s_{m(i,j)}(\cdot ;X_{i-1,j}^{\mathbf
x},X_{i,j-1}^{\mathbf x})$ to depend on $U_{l,n}$, $(l,n) \in
R_{h,k}$, with $l\leq i$ and $n \leq j$, in such a way that
$$
s_{m(i,j)}(U_{l,n},\,\,(l,n) \in R_{h,k}, l\leq i, n \leq
j;X_{i-1,j}^{\mathbf x},X_{i,j-1}^{\mathbf x} )\sim
K(\cdot|X_{i-1,j}^{\mathbf x},X_{i,j-1}^{\mathbf x}).
$$
A particular example is when $s_{m(i,j)}$ realizes the maximal
coupling of the laws $\{K(\cdot|x_1,x_2)$, $(x_1,x_2) \in
I_{i-1,j}\times I_{i,j-1}\}$ where $I_{i,j}$ is defined as in
(\ref{range}).

Now consider any finite subset $S \subset {\mathbb Z}^2$. By
analogy with Lemma \ref{R2} it is quite clear that if the set $B
\subset S$ is such that $\partial_I B$ is contained in a union of
blocks $B_{h,k}$ where $W_{h,k}=0$, then the restriction of the
field ${\mathbf X}^{\mathbf x}_S$ on $B$ can be represented as a
measurable function $\Upsilon^B({\mathbf U}_B)$, irrespectively of
the b.c.'s ${\mathbf x} \in {\overrightarrow{\partial} S}$.

We say that $\gamma=\{(i_k, j_k) \in {\mathcal V}_d\}_{k=0,1,
\dots,m}$ is an increasing path of length $m$ if for
$k=0,\ldots,m-1$ the vertex $(i_k, j_k)$ is a parent of $(i_{k+1},
j_{k+1})$ in the graph ${\mathcal G}_d$. The field ${\mathbf W}$
is then used to define open paths in ${\mathcal G}_d$.

Notice that a field ${W}_{h,k}=\psi (\mathbf U_{R_{h,k}})$ defined
for $(h,k) \in {\mathcal V}_d$ (see Assumption \ref{assume3}) is
not a Bernoulli field in general. However since $\mathbf U$ is an
i.i.d. field the random variables
$$
\{W_{h_i,k_i}, (h_i,k_i) \in {\mathcal V}_d, \,\,i=1,\ldots, t\}
$$
are mutually independent whenever the regions $R_{h_i,k_i}$ are
pairwise disjoint, for $i=1,\ldots,t$. Now it can be verified that
$R_{h_1,k_1}$ and $R_{h_2,k_2}$ are not disjoint if and only if
$h_1+k_1=h_2+k_2=j(d-1)$ for some integer $j$ and $|h_1-h_2|\leq
d-1$ (e.g. check in Figure \ref{fig1} that $R_{0,0} \bigcap
R_{1,-1} = \{(2,0)\}$). Therefore the field on each diagonal is
$(d-1)$-dependent and different diagonals are mutually
independent.

For the sake of simplicity, in this section we suppose that
$\partial_I {\Lambda}$ is a union of blocks $B_{h,k}$, so it can
be identified with a subset of ${\mathcal V}_d$; otherwise we can
enlarge $\Lambda$ to have this property. In analogy with
(\ref{verticic}) we define the random region
\begin{equation}\label{vert}
\sigma(\Lambda) = \sigma(\Lambda,{\mathbf W})= \{ (i,j)\in
{\mathcal V}_d: \exists \hbox{ an open path in ${\mathcal G}_d$
joining } (i,j) \hbox{ to } (k,l) \in
\partial_{I} \Lambda \}.
\end{equation}
and $D(\Lambda)=\Lambda \cup \sigma(\Lambda)$. As in Lemma
\ref{uscita} we can prove that for any $(h,k) \in \partial_I
D(\Lambda)$ the random variable $W_{h,k}=0$. We can finally show
the following

\begin{theorem}\label{generaliz}
Suppose the kernel $K$ satisfies Assumption \ref{assume3}. Then
$K$ is uniformly ergodic and  for any finite region $\Lambda$ the
law of its stationary version ${\mathbf X'}_{\Lambda}$ is the same
as $\Upsilon^{D(\Lambda)}({\mathbf U}_{D(\Lambda)})$. Moreover
${\mathbb E} |D (\Lambda )| =| \Lambda | + O ( | \partial_I
\Lambda | )$.
\end{theorem}
\begin{proof} As in (\ref{verticic})
$$
\sigma(\Lambda)=\bigcup_{(h,k) \in \partial_I \Lambda} C'_{h,k}
$$
where $C'_{h,k}$ is the set of blocks which can be joined to
$(h,k)$ by an open path in ${\mathcal G}_d$. We have only to prove
that under Assumption \ref{assume3} $C'_{0,0}$ is finite almost
surely and has finite mean. The only difference with Theorem
\ref{main} is that ${\mathbf W}$ is not Bernoulli. However, since
different diagonals are independent, the values of the field
${\mathbf W}$ on any increasing path are i.i.d. Next define
\begin{equation}\label{diametr}
  \hbox{diam}(C'_{0,0} ) = \max \{
  \hbox{ length of an open path from $(h,k)$ to $(0,0)$: } (h,k ) \in C'_{0,0} \} .
\end{equation}
Then the following holds
\begin{equation}\label{pata}
  P( |C'_{0,0}| \geq k^2  )\leq P ( \hbox{diam}(C'_{0,0} )\geq k )
  \leq \{d (1- {\bar \delta })\}^k ,
\end{equation}
where the first inequality is due to a trivial geometric argument.
The second inequality follows since the probability that a fixed
increasing path of length $k $ is open is equal to $ (1- {\bar
\delta })^k$ and the number of such paths leading to the origin is
equal to $d^k$. Since by Assumption \ref{assume3} $d (1- {\bar
\delta }) <1 $ we obtain that $ C'_{0,0}$ is finite w.p. $1$ and
its cardinality has finite mean (in fact, finite moments of any
order).
\end{proof}

It is natural to ask whether it is possible to improve the
inequality (\ref{jam}) by dominating $\mathbf W$ with a suitable
Bernoulli field. A domination result of this type for
$(d-1)$-dependent fields can be found in \cite{LSS}: however we
have checked that it does not improve the bound (\ref{jam}).

Before discussing some examples we need to make clear the
limitations of the block algorithm presented in this section, due
to the difficulties arising from the two-dimensional structure. In
some sense we have constructed a skeleton process $\mathbf
X^{\mathbf x}_{B_{h,k}}$, $(h,k) \in {\mathcal V}_d$: however,
this field is not a ${\mathcal G}_d$-unilateral field anymore, due
to the overlap between the regions $R_{h,k}$, for $(h,k) \in
{\mathcal V}_d$. As a consequence a minorization condition of the
type
$$
P(\mathbf X^{\mathbf x}_{B_{h,k}}\in A|\mathbf
  \mathbf X^{\mathbf x}_{\overrightarrow{\partial} R_{h,k}}={\mathbf x})\geq \delta \varphi(A)
$$
for all ${\mathbf x} \in E^{\overrightarrow{\partial} R_{h,k}}$,
does not  immediately translates in the possibility of coupling
the field with an auxiliary Bernoulli field with probability
$\delta$ of $0$. For this reason we need the stronger Assumption
\ref{assume3}, which requires the choice of a specific family of
functions $\{s_m\}$ allowing to simulate the field also at sites
common to different regions. Moreover, since the binary auxiliary
field used by the algorithm is not Bernoulli, the values of
$\delta$ for which we can prove that the algorithm works have to
be larger than those obtained with i.i.d. percolation ($\delta
> 1/2$ rather than the previously cited bound
$ \delta \geq 0.317 \dots $, for $d =2$).

As a consequence we cannot identify the class of uniformly ergodic
fields which satisfy Assumption~\ref{assume3} for some choice of
$l$ and $d$. As a matter of fact we are not aware of any simple
characterizations of the class of uniformly ergodic unilateral
random fields, as possible in the one-dimensional case.

\section{Examples} \label{aria}

\noindent {\bf Example 1.} In this example we will construct a
perfect simulation algorithm for the stationary unilateral field
with kernel $K$ in a finite region under the following

\begin{assumption}\label{assume2}
There exists a measurable subset $C \subset E$, a probability
measure $\phi$ on $E$ and positive constants $0<\rho_1<1 $ and
$0<\rho_2<1$ with $\rho_1^2 \rho_2 > 2/3$ such that
\begin{equation}\label{pric}
K(C |y_1, y_2) \geq \rho_1 \,\,\,\,\,\,\, \forall \,  (y_1, y_2)
\in E^2
\end{equation}
and
\begin{equation}  \label{ddoblin}
K(A|y_1,y_2)\geq \rho_2 \phi (A)\,\,\,\,\,\,\, \forall \,  (y_1,
y_2) \in C^2
\end{equation}
for every measurable set $A$.
\end{assumption}

Under Assumption \ref{assume2} we will build  a coupling of the
field with kernel $K$ which satisfies Assumption \ref{assume3}
with $l=1$ and $d=3$. In this particular case we need to
distinguish between even and odd numbered diagonals, depending on
the parity of the sum of the coordinates of a site: the blocks
will be single vertices $(h,k)$ lying on even numbered diagonals.

The field will be constructed by means of two families of
functions $s_0$ and $s_1$ with the property (\ref{couple}), which
are used to get the value at sites lying on even and odd numbered
diagonals, respectively. For simplicity of notation we use pairs
of random variables $(U_1,U_2)$ which are uniformly distributed in
$(0,1)^2$ rather than a single random variable in $(0,1)$. The
functions $s_0$ and $s_1$ are defined by means of some functions
$f_i:(0,1) \times E^2\rightarrow E$, for $i=0,1,2,3,4$, measurable
in each of the arguments, with the following properties. For
$(y_1,y_2)\in E^2$
\begin{equation}
f_0(U_1;y_1,y_2)=f_0(U_1) \sim \phi(\cdot),
\end{equation}
\begin{equation}
f_1 ( U_1; y_1, y_2 ) \sim \left \{
\begin{array}{c}
 \frac{1}{1-\rho_2} [ {K(\cdot|y_1,y_2)-\rho_2
\phi(\cdot)}] \hbox{ for } y_1,y_2 \in C  \\
K(\cdot|y_1,y_2) \hbox{ otherwise.}
\end{array}    \right .
\end{equation}
\begin{equation}
f_2 ( U_1; y_1, y_2 ) \sim {K ( \cdot | y_1, y_2 )}.
\end{equation}
\begin{equation}
f_3 ( U_1; y_1, y_2 ) \sim \frac{K ( \cdot\cap C | y_1, y_2 )}{ K
( C| y_1, y_2 )} ,
\end{equation}
\begin{equation} \label{46a}
f_4 ( U_1; y_1, y_2 ) \sim \frac{1}{1-\rho_1}\left \{ K( \cdot |
y_1, y_2 ) - \rho_1 \frac{K ( \cdot\cap C | y_1, y_2 )}{ K ( C|
y_1, y_2 )} \right \}.
\end{equation}
These functions always exist since $E$ is assumed to be Borel.
Notice that (\ref{46a}) is well defined since
\begin{equation}\label{oaoa}
   K( \cdot |
y_1, y_2 ) - \rho_1 \frac{K ( \cdot\cap C | y_1, y_2 )}{ K ( C|
y_1, y_2 )}  = K ( \cdot\cap C^c | y_1, y_2 ) +\left ( 1
-\frac{\rho_1}{K (C | y_1, y_2 )} \right ) K ( \cdot\cap C | y_1,
y_2 )
\end{equation}
is non-negative by (\ref{pric}).

With these positions it is not difficult to verify that the
functions
\begin{equation}\label{fung0}
  s_0 ( u_1,u_2 ; y_1, y_2) = \left \{ \begin{array}{ll}
     &(1-1_{(0,\rho_2)}(u_2))  f_0(u_1)
     +1_{(0,\rho_2)}(u_2)   f_1  (u_1; y_1, y_2)
    \hbox{ for }  y_1, y_2 \in C,\\
  &f_2  (u_1; y_1, y_2
  )     \hbox{ otherwise. }
  \end{array}
   \right .
\end{equation}
\begin{equation}\label{fung1}
  s_1 ( z, u_1,u_2 ; y_1, y_2) = 1_{(0,\rho_1)}(u_2) f_3(u_1; y_1, y_2 )
  +(1-1_{(0,\rho_1)}(u_2)) f_4  (u_1; y_1, y_2) ,
\end{equation}
satisfy (\ref{couple}). Essentially $s_0$ and $s_1$ come from two
different mixture decompositions of the kernel. Finally the field
${\mathbf W}$ is defined for any pair $(h,k) \in {\mathbb Z}^2$
such that $(h+k)$ is even by
\begin{equation}\label{Wf}
  W_{h,k} =
    1- 1_{(0,\rho_1)}((U_2)_{h-1,k})1_{(0,\rho_1)}((U_2)_{h,k-1})
    1_{(0,\rho_2)}((U_2)_{h,k}).
\end{equation}
Notice that $W_{h,k}=0$ if and only if all the indicators
appearing in the above formula are equal to $1$, which implies
that $X_{h,k}=f_0((U_1)_{h,k})$, irrespectively of the values at
the parent sites. Since
$$
P(W_{i,j}=0)=\rho_1^2 \rho_2>\frac{2}{3}
$$
we have thus established
\begin{theorem}\label{T2}
  Under Assumption~\ref{assume2} the conclusions of Theorem \ref{generaliz}
  hold.
\end{theorem}

It is not difficult to construct a kernel $K$ for which Assumption
\ref{assume2} is satisfied, whereas Assumption \ref{assume} is
false. Let $E=\{0,1,2\}$ and suppose that
\begin{equation} \label{condizq}
K(l|y_1,y_2)= \left \{ \begin{array}{l}
   \phi (l) \hbox{ if }  (y_1 , y_2) \in \{0,1\}^2\\
  \delta_{1,l} \hbox{ if } \,\,\,\, y_1=y_2=2 \\
\delta_{0,l} \,\,\,\, \hbox{ otherwise. }
\end{array}  \right .
\end{equation}
Choosing $C=\{0,1\}$ (\ref{pric}) holds with
$\rho_1=\phi(C)=\phi(0)+\phi(1)$ and (\ref{ddoblin}) holds with
$\rho_2=1$, hence if $\phi(C)$ is large enough Assumption
\ref{assume2} is satisfied. On the other hand Assumption
\ref{assume} does not hold, since $K(\cdot|2,2)$ and
$K(\cdot|1,2)$ are singular.

\noindent {\bf Example 2.} This example allows a simple analysis
of the role of $d$ and $l$ in the block algorithm. On the state
space $E=\{0,1,2\}$ consider the unilateral kernel
\begin{equation}\label{kenel}
  K(l|y_1,y_2)=\left \{ \begin{array}{l}p\delta_{2,l}+(1-p)\delta_{1,l}
   \hbox{ if } y_1=y_2=2 \\
    p\delta_{1,l}+(1-p)\delta_{0,l}
   \hbox{ if } \min(y_1,y_2)=0\\
    p\delta_{\min(y_1,y_2)+1,l}+(1-p)\delta_{\min(y_1,y_2)-1,l} \hbox{ otherwise. }
\end{array} \right .
\end{equation}
It is clear that this kernel does not satisfy Assumption
\ref{assume} for any $p \in [0,1]$, since
$K(0|2,2)=K(1|1,1)=K(2|0,0)=0$. Now assume $p<1/2$ and consider
what happens by increasing $d$ or $l$. For $d=3$ and $l=1$ we see
that by choosing $C=\{0,1\}$ and $\phi$ the Dirac mass on $0$,
with $\rho_1=\rho_2=1-p$, we can fulfil Assumption \ref{assume2}
provided $(1-p)^3>2/3$. This choice of $C$ is clearly the best
possible.

Next we show that by taking $l=2$ and $d=2$ we can enlarge the
region of parameters $p$ under which uniform ergodicity can be
proved and a perfect simulation algorithm can be constructed
through Assumption \ref{assume3}. For any site of the lattice we
always use the natural coupling $s:(0,1)\times E^2 \to E$, defined
by
\begin{equation}\label{trivia}
  s(u;y_1,y_2)=1_{(0,p)}(u)\min\{2,\min(y_1,y_2)+1\}+1_{(p,1)}(u)\max\{0,\min(y_1,y_2)-1\}.
\end{equation}
We define $W_{0,0}=0$ provided $U_{i,j}>p$ on all the sites
$(i,j)=(1,2),(1,1),(2,1)$. The whole field $\{W_{h,k}, (h,k) \in
{\mathbb Z}^2\}$ is defined by translation. Such a definition
ensures that the field has the value zero on both the sites
$(2,1)$ and $(1,2)$ which form the block $B_{0,0}$, irrespectively
of the b.c.'s on the parent blocks, formed by the sites
$(-1,2),(0,1),(1,0),(2,-1)$ (see again Figure \ref{fig1}). Hence
we can take $\Phi=0$ in (\ref{ham}). It is immediately obtained
that $P(W_{0,0}=0)=(1-p)^3$, so that Assumption \ref{assume3} is
satisfied provided $(1-p)^3>1/2$, enlarging the region where we it
is proved that perfect simulation works.


\begin{thebibliography}{9}

\bibitem{Du84}  R. Durrett, Oriented percolation in two dimensions. {\it
Ann. Probab.} {\bf 12}: 999-1040 (1984).

\bibitem{Ferra}
P.A. Ferrari, R. Fern\'andez, N.L. Garcia, Perfect simulation for
interacting point processes, loss networks and Ising models. {\it
Stoch. Proc. Appl.} {\bf 102}: 63--88 (2002).

\bibitem{FT} S. Foss, R.L. Tweedie, Communication in statistics,
stochastic models, {\it Comm. Statist. Stoch. Models} {\bf 14}:
187--203 (1998).

\bibitem{GW1} R.F. Galbraith, D. Walley, Ergodic properties of a
two-dimensional binary processes, {\it J. Appl. Prob.} {\bf 17}:
124--133 (1980).

\bibitem{GW2} R.F. Galbraith, D. Walley, Further properties for
unilateral binary processes, {\it J. Appl. Prob.} {\bf 19}:
332--343 (1982).

\bibitem{G}
J. Goutsias, Unilateral approximation of Gibbs random field
images, {\it Computer Vision, Graphics, and Image Processing:
Graphical Models and Image Processing} {\bf 53}: 240--257, (1991).

\bibitem{Gra}
L. Gray, J.C. Wierman, R.T. Smythe, Lower bounds for the critical
probability in percolation models with oriented bonds, {\it J.
Appl. Probab.} {\bf 17}: 979--986, (1980).


\bibitem{Hag}
O. H\"aggstr\"om, Finite Markov chains and algorithmic
applications. Mathematical Society Student Texts, {\bf 52}.
Cambridge University Press, Cambridge, 2002.

\bibitem{Ki}
Y. Kifer, Ergodic theory of random transformations,
 Progress
   in Probability and Statistics, {\bf 10}.  {
Birkh\"auser}, Boston, 1986.

\bibitem{L}
S.L. Lauritzen, Graphical models. Oxford Statistical
   Science Series, {\bf 17}. Oxford Science Publications, New York,
   1996.

\bibitem{LSS}
T.M. Liggett, R.H. Schonmann, A.M. Stacey, Domination by product
measures. {\it Ann. Probab.} {\bf 25}: 71-95 (1997).

\bibitem{Mad}
N. Madras (ed.), Monte Carlo methods. Fields Institute
Communications,  {\bf 26}. AMS, Providence, 2000.

\bibitem{MG}
D.J. Murdoch, P.J. Green,
 Exact sampling from a continuous state space. {\it Scand.
   J. Statist.} {\bf 25}:  483--502 (1998).

\bibitem{MT}
S.P. Meyn, R.L. Tweedie, Markov chains and stochastic stability.
Communications and Control Engineering Series. {\it
Springer-Verlag}, London, 1993.

\bibitem{Pi1}
D.K. Pickard, Unilateral Markov fields. {\it Adv. Appl. Prob.}
{\bf 12}: 655--671 (1980)

\bibitem{PW}
J.G. Propp, D.B. Wilson, Exact sampling
   with coupled Markov chains and applications to statistical
   mechanics.
   {\it Random Structures Algorithms} {\bf 9}: 223--252 (1996)

\bibitem{WG73}
T.R. Welberry, R. Galbraith, A two-dimensional model of crystal
growth. {\it J. Appl. Cryst.} {\bf 6}: 87--96 (1973).

\end{thebibliography}
\end{document}